\documentclass[11pt]{article}%

\usepackage[utf8]{inputenc}

\usepackage{amsmath}

\usepackage{amsfonts}
\usepackage{mathrsfs}
\usepackage{amssymb, color}
\usepackage[linkcolor=black,anchorcolor=black,citecolor=black]{hyperref}
\usepackage{graphicx}
\usepackage{accents}
\numberwithin{equation}{section}
\usepackage[body={15.5cm,21cm}, top=3cm]{geometry}%
\providecommand{\U}[1]{\protect\rule{.1in}{.1in}}
\providecommand{\U}[1]{\protect \rule{.1in}{.1in}}
\newtheorem{theorem}{Theorem}[section]

\newtheorem{example}[theorem]{Example}

\newtheorem{lemma}[theorem]{Lemma}

\newtheorem{assumption}[theorem]{Assumption}

\begin{document}
	\title{Stochastic maximum principle for optimal control problem with varying terminal time and non-convex control domain}
	\author{Jin Shi\thanks{Shandong University-Research Center for Mathematics  and Interdisciplinary Sciences, Shandong University, PR China,
(202011968@mail.sdu.edu.cn).} \quad  Shuzhen Yang \thanks{Shandong University-Zhong Tai Securities Institute for Financial Studies, Shandong University, PR China, (yangsz@sdu.edu.cn). }}
	\date{}
	\maketitle

\textbf{Abstract}: In this paper, we consider a varying terminal time structure for the stochastic optimal control problem under state constraints, in which the terminal time varies with the mean value of the state. In this new stochastic optimal control system, the control domain does not need to be convex and the diffusion coefficient contains the control variable. To overcome the difficulty in the proof of the related Pontryagin’s stochastic maximum principle, we develop asymptotic first- and second-order adjoint equations for the varying terminal time, and then establish its variational equation. In the end, two examples are given to verify the main results of this study.
	
		
{\textbf{MSC2010}: 93E03; 93E20; 60G99
		
\addcontentsline{toc}{section}{\hspace*{1.8em}Abstract}
	
\section{Introduction}

In the classical stochastic optimal control problem, we would like to minimize the following cost functional,
\begin{equation}
\label{incos-1}
J(u(\cdot))=\mathbb{E}\bigg{[}\displaystyle\int_0^Tf(X^u(t),u(t))\mathrm{d}t+g(X^u(T))\bigg{]},
\end{equation}
where $f\left(X(t),u(t)\right)$ is the running cost at time $t\in \left[0,T\right]$ with a given terminal time $T$, $g(X^u(T))$ is the terminal cost, and the state process $X^{u}\left(t\right)$ satisfies
\begin{equation}
\label{insde-1}
X^u(s)=x_0+\int_{0}^{s}b(X^u(t),u(t))\mathrm{d}t+\int_{0}^{s}\sigma(X^u(t),u(t))\mathrm{d}W(t),\ s\in [0,T].
\end{equation}
To study the above optimal controlled problem, stochastic maximum principle is developed to establish the related necessary condition: when diffusion coefficient $\sigma$ contains the control variable $u$ and the control domain is convex, to obtain the local stochastic maximum principle, the convex variation method was introduced in Besoussan\cite{B81} and Bismut \cite{B78}; when the control domain does not convex, and diffusion coefficient contains the control variable, Peng \cite{P90} first introduced the second-order adjoint equation and established the stochastic maximum principle in global form by using the "spike variation method". \cite{P90} overcame the difficulty in the stochastic optimal control problem where the diffusion coefficient contains the control variable and the control domain is non-convex. Since then, many authors begin to consider the global stochastic maximum principle of some general cases: the global stochastic maximum principle for the stochastic recursive utilities system was established in Hu \cite{17H}; the stochastic maximum principle for fully coupled forward-backward stochastic system was considered in Hu and Ji \cite{16HJ}. Also see \cite{BE10,FM06,F10,W13,Y18,Y10,Y99} for the theory of stochastic optimal control problems.

In this paper, we focus on the stochastic optimal control problem under state constraints, where a varying terminal time optimal control structure is given by
\begin{equation}
\label{intime-1}
\tau^u=\inf\bigg{\{}t:\mathbb{E}[\Phi(X^u(t))]\geq \alpha,\ t\in [0,T] \bigg{\}}\bigwedge T.
\end{equation}
We would like to minimize the cost functional
\begin{equation}
\label{incos-2}
J(u(\cdot))=	\mathbb{E}\bigg{[}\displaystyle\int_{0}^{\tau^u}f(X^u(t),u(t))\mathrm{d}t+g(X^u(\tau^{u}))\bigg{]},
\end{equation}
subject to
\begin{equation}
\label{insde-2}
\text{d}{X}^u(t)=b(X^u{(t)},u(t))\text{d}t+\sigma (X^u{(t)},u(t))\text{d}W(t).
\end{equation}

For this novel stochastic optimal control problem with a varying terminal time, Yang \cite{Y20} first established the local stochastic maximum principle, and the method used in \cite{Y20} is the convex variation method. Then, the varying terminal time mean-variance model was solved in Yang \cite{Y22}. Furthermore, Wang and Yang \cite{JS24} developed the stochastic maximum principle for the recursive optimal control problem with varying terminal time. For the optimal control problem with a varying terminal time, when the control domain is non-convex, the related global stochastic maximum principle is still unsolved.

In this study, we want to establish the related global stochastic maximum principle for the optimal control problem with a varying terminal time. Following the idea given in Peng \cite{P90} and Yang \cite{Y20}, we will use the spike variation method to develop the stochastic maximum principle. Precisely, let $\varepsilon>0$, $\bar{u}(\cdot)$ be an optimal control and $u\in U$, where $U$ is the control domain. We consider the following control:
\begin{equation}
u^{\varepsilon}\left(t\right)=\begin{cases}
u, &\text{if} \quad  \tau\leq t \leq \tau+\varepsilon ;\\
\bar{u}\left(t\right), & \text{otherwise}.
\end{cases}
\end{equation}
This $u^{\varepsilon}$ is called a spike variation of the optimal control $\bar{u}(\cdot)$. We denote the optimal  pair for cost functional$(\ref{incos-2})$  by $(\bar{u}(\cdot),\bar{X}(\cdot))$.
To overcome the difficulty in the proof of the stochastic maximum principle for the stochastic optimal control problem, we develop asymptotic first- and second-order adjoint equations for the varying terminal time, and establish the related variational equations for varying terminal time and cost functional based on three kinds of optimal terminal time $\tau^{\bar{u}}$: $\text{(i)} \ \tau^{\bar{u}}<T$;
$\text{(ii)}\ \inf\bigg{\{}t:\mathbb{E}[\Phi(\bar{X}(t))]\geq \alpha,\ t\in [0,T] \bigg{\}}=T$;
$\text{(iii)}\ \bigg{\{}t:\mathbb{E}[\Phi(\bar{X}(t))]\geq \alpha,\ t\in [0,T] \bigg{\}}=\varnothing$, and then obtain the global stochastic maximum principle.

The remainder of this paper is organized as follows: in Section 2, we describe a stochastic optimal control problem with a varying terminal time. Then, we give some preliminary results for proving the global stochastic maximum principle in Section 3, and we establish the global stochastic maximum principle in Section 4. We conclude the main results in Section 5.

\section{The optimal control problem}

In this section, we introduce the stochastic control problem with a varying terminal time and non-convex control domain. Given a complete filtered probability space $(\Omega,\mathcal{F},P;\{ \mathcal{F}(t)\}_{t\geq0})$, where  the  filtration $\{ \mathcal{F}(t)\}_{t\geq0}$ is generated by $d$-dimensional standard Brownian motion $W$. Let $T>0$ be a given constant, control domain $U$ be a nonempty subset of $\mathbb{R}^k$  with a given positive integer $k$ and the set of all admissible controls $\mathcal{U}[0,T]=\left\{u(\cdot)\in\mathcal{A}^{2}|u(t)\in U, 0\leq t\leq T,\ a.s.\right\}$, where $\mathcal{A}^{2}$ denotes the set of all square integrable and  $\mathcal{F}(t)$-adapted processes.

We introduce the following controlled stochastic differential equation with $u(t)\in\mathcal{U}[0,T]$,
\begin{equation}
\label{insde-3}
\text{d}{X}^{u}(t)=b(X^{u}(t),u(t))\text{d}t+\sigma(X^{u}(t),u(t))\text{d}W(t) ,\quad t\in(0,T],
\end{equation}
and $X(0)=x_0$, where
\[%
\begin{array}
	[c]{l}%
	b:\mathbb{R}^m\times U\to \mathbb{R}^m,\\
	\sigma:\mathbb{R}^{m}\times U\to \mathbb{R}^{m\times d}.\\
\end{array}
\]

In this study, we minimize the following varying terminal time cost functional over admissible  control set $\mathcal{U}[0,\tau^{u}]$:
\begin{equation}	
\begin{aligned}
\label{incos-3}
&J(u(\cdot))=\mathbb{E}\bigg{[}{\displaystyle \int \limits_{0}^{\tau^{u}}}
f(X^u{(t)},u(t))\text{d}t+g(X^u(\tau^u))\bigg{]},
\end{aligned}
\end{equation}
where
\begin{equation}
\label{time-2}
\tau^u=\inf\bigg{\{}t:\mathbb{E}[\Phi(X^u(t))]\geq \alpha,\ t\in [0,T] \bigg{\}}\bigwedge T,
\end{equation}
and
\[%
\begin{array}
	[c]{l}%
f:\mathbb{R}^m\times U\to \mathbb{R},\\
\Phi,g:\mathbb{R}^{m}\to \mathbb{R},\\
\end{array}
\]
$a_1\bigwedge a_2=\min(a_1,a_2),\ a_1,a_2\in \mathbb{R}$ and $\alpha\in (\Phi(x_0),+\infty)$. Note that if $\alpha<\Phi(x_0)$, then $\tau^u=0$, and the problem is trivial.

We set $\sigma=(\sigma^1,\sigma^2,\cdots,\sigma^d)$, and $\sigma^j\in \mathbb{R}^m$ for $j=1,2,\cdots, d$. Here,
$\mathbb{R}^m=\mathbb{R}^{m\times 1}$. In addition, "$*$" denotes the transform of vector or matrix.

To investigate the stochastic optimal control problem, we introduce the following assumptions.
\begin{assumption}
\label{ass-a}
Let $b$, $\sigma$, $f$, $g$ be twice continuously differentiable with respect to $x$. $b$, $\sigma$, $f$ are continuous with respect to $u$. They and all their derivatives $b_{x}$, $b_{xx}$, $\sigma_{x}$, $\sigma_{xx}$, $f_{x}$, $f_{xx}$, $g_{x}$, $g_{xx}$ are continuous in $\left(x,u\right)$. $b_{x}$, $b_{xx}$, $\sigma_{x}$, $\sigma_{xx}$, $f_{xx}$, $g_{xx}$ are bounded, and $b$, $\sigma$, $f_{x}$, $g_{x}$ are bounded by $C\left(1+\left | x \right |+\left | u \right |  \right)$.
\end{assumption}

\begin{assumption}
\label{ass-b}
$\Phi$ are twice continuously differentiable with respect to $x$.
\end{assumption}

 The correspond solution $X^{u}\left(\cdot\right)$ is called state variable or trajectory under control $u\left(\cdot\right)$, and the optimal control $\bar{u}(\cdot)$ satisfying
\begin{equation}
J(\bar{u}(\cdot))= \underset{u(\cdot)\in\mathcal{U}[0,\tau^u]}{\inf}J(u(\cdot)).
\end{equation}
The corresponding state trajectory $(\bar{u}(\cdot),\bar{X}(\cdot))$ is called an optimal state trajectory or optimal pair and $\tau^{\bar{u}}$ is  called the optimal terminal time.

\section{Preliminary results}
To obtain the well-known Pontryagin's stochastic maximum principle for the varying terminal time optimal control problem, we show the preliminary results in this section. Due to the appearance of the control variable in $\sigma\left(\cdot,\cdot\right)$ and the non-convex control domain $U$, the usual first-order expansion approach does not work. Hence, we need to introduce a second-order expansion method. The second-order expansion method was first proposed by Peng in \cite{P90}. We construct an admissible control in the following way:
\begin{equation}
u^{\varepsilon}\left(t\right)=\begin{cases}
u, & \text{if} \ \tau \leq t \leq \tau+\varepsilon ,\\
\bar{u}\left(t\right), & \text{otherwise},
\end{cases}
\end{equation}
where $0 \leq \tau <T$ is fixed, and $\varepsilon>0$ is sufficiently small. $u$ is an arbitrary $\mathcal{F}_{\tau}$-measurable bounded random variable with values in $U$. Let $X^{\varepsilon}\left(\cdot\right)$ be the trajectory of the control system $\label{insde-3}$ corresponding to the control $u^{\varepsilon}\left(\cdot\right)$.

For the function $\Phi(X^u(\cdot))$ in constrained condition (\ref{time-2}), we have the following results which was given in \cite{Y20}.
\begin{lemma}
\label{le-1}
Let Assumptions \ref{ass-a} and \ref{ass-b} hold. We have
\begin{equation}
\label{mean-0}
\mathbb{E}[\Phi(X^u(s))]=\Phi(x_0)+\int_0^sh^u(t)\mathrm{d}t,\ s\in [0,T],
\end{equation}
where
$h^{{u}}(t)=\mathbb{E}\bigg{[}\Phi_x(X^{{u}}(t))^{*}b(X^{{u}}(t),{u}(t))+\displaystyle\frac{1}{2}\sum_{j=1}^d\sigma^j(X^{{u}}(t),{u}(t))^{*}\Phi_{xx}(X^{{u}}(t))\sigma^j(X^{{u}}(t),{u}(t))\bigg{]}$, $t\in [0,s]$.	
\end{lemma}

For the optimal terminal time, we have the following asymptotic behavior. The proof of this lemma is similar to Lemma 2 in Yang \cite{Y20}. Thus, we omit the proof.
\begin{lemma}
\label{le-2}
Let Assumptions \ref{ass-a} and \ref{ass-b} hold, and suppose that $h^{\bar{u}}(\tau^{\bar{u}})\neq 0$ and $h^{\bar{u}}(\cdot)$ is continuous at the point $\tau^{\bar{u}}$. We have the following results.
	
 (i). If $\tau^{\bar{u}}<T$, one obtains
\begin{equation}
\label{le0-e0}
\lim_{\varepsilon\to 0}\left|{\tau^{\bar{u}}-\tau^{u^{\varepsilon}}}\right|=0.
\end{equation}

(ii). If $\inf\bigg{\{}t:\mathbb{E}[\Phi(\bar{X}(t))]\geq \alpha,\ t\in [0,T] \bigg{\}}=T$, one obtains
\begin{equation}
\lim_{\varepsilon\to 0}\left|{\tau^{\bar{u}}-\tau^{u^{\varepsilon}}}\right|=0.
\end{equation}

(iii). If $\bigg{\{}t:\mathbb{E}[\Phi(\bar{X}(t))]\geq \alpha,\ t\in [0,T] \bigg{\}}=\varnothing$, we have
\begin{equation}
\lim_{\varepsilon\to 0}\left|{\tau^{\bar{u}}-\tau^{u^{\varepsilon}}}\right|=0.
\end{equation}	
\end{lemma}

For the optimal pair, we have the following asymptotic behavior which is given in \cite{P90}, we omit the proof.
\begin{lemma}
\label{le-3}
Let Assumptions \ref{ass-a} and \ref{ass-b} hold, we have
\begin{equation}
		\varepsilon ^{-2}\underset{0\leq t\leq T}\sup\mathbb{E}\left | X^{\varepsilon}\left (t  \right )-\bar{X}(t)-y_{1}\left ( t \right ) -y_{2}\left ( t \right )  \right |^{2} \le  C,
\end{equation}
where $y_{1}\left(\cdot\right)$,\ $y_{2}\left(\cdot\right)$ are solutions of the following the first- and second-order variational equations:
\begin{equation}
\label{eq-1}
\begin{aligned}
			y_{1}\left(t\right)=& \int_{0}^{t} b_{x}\left (\bar{X}\left ( s \right ),\bar{u}\left ( s \right )   \right )y_{1}\left(s\right)+\left(b\left (\bar{X}\left ( s \right )  ,u^{\varepsilon}\left ( s \right )  \right )-b\left (\bar{X}\left ( s \right ),\bar{u}\left ( s \right )   \right )\right) \mathrm{d}s\\
			&+\int_{0}^{t} \sigma _{x}\left (\bar{X}\left ( s \right ),\bar{u}\left ( s \right )   \right )y_{1}\left(s\right)+\left(\sigma \left (\bar{X}\left ( s \right )  ,u^{\varepsilon}\left ( s \right )  \right )-\sigma\left (\bar{X}\left ( s \right ),\bar{u}\left ( s \right )   \right )\right) \mathrm{d}W\left ( s \right ),
\end{aligned}
\end{equation}

\begin{equation}
\label{eq-2}
\begin{aligned}
		    y_{2}\left(t\right)=& \int_{0}^{t} b_{x}\left (\bar{X}\left ( s \right ),\bar{u}\left ( s \right ) \right )y_{2}\left(s\right)+\frac{1}{2}b_{xx}\left (\bar{X}\left ( s \right ),\bar{u}\left ( s \right ) \right )y_{1}\left(s\right)y_{1}\left(s\right)\mathrm{d}s\\	
			&+\int_{0}^{t} \sigma _{x}\left (\bar{X}\left ( s \right ),\bar{u}\left ( s \right )   \right )y_{2}\left(s\right)+\frac{1}{2}\sigma_{xx}\left (\bar{X}\left ( s \right ),\bar{u}\left ( s \right ) \right )y_{1}\left(s\right)y_{1}\left(s\right)\mathrm{d}W\left(s\right)\\	
			&+\int_{0}^{t}\left(b_{x}\left (\bar{X}\left ( s \right )  ,u^{\varepsilon}\left ( s \right )  \right )-b_{x}\left (\bar{X}\left ( s \right ),\bar{u}\left ( s \right )  \right )\right)y_{1}\left(s\right)\mathrm{d}s\\
			&+\int_{0}^{t}\left(\sigma_{x} \left (\bar{X}\left ( s \right )  ,u^{\varepsilon}\left ( s \right )  \right )-\sigma_{x}\left (\bar{X}\left ( s \right ),\bar{u}\left ( s \right )   \right )\right)y_{1}\left(s\right) \mathrm{d}W\left ( s \right ),\\	
\end{aligned}
\end{equation}
where $r_{xx}yy=\sum_{i,j=1}^{m}r_{x^{i}x^{j}}y^{i}y^{j} $, for $r=b,\sigma, f,g.$
\end{lemma}

To investigate the stochastic maximum principle under a varying terminal time for non-convex control domain, we need to develop the variational equation for optimal terminal time.
For notation simplicity, we denote $\bar{h}(u,t):=\lim_{\varepsilon\to 0}\int_{0}^{\tau^{{u}^{\varepsilon}}}\frac{h^{u^{\varepsilon}}\left(t\right)-h^{\bar{u}}\left(t\right)}{\varepsilon}dt.$
We show the existence of $\bar{h}(u,t)$ in the following lemma.

Denote by
$$
l\left(X^{{u}}(t),{u}(t)\right):=\Phi_x(X^{{u}}(t))^{\top}b(X^{{u}}(t),{u}(t))
+\displaystyle\frac{1}{2}\sum_{j=1}^d\sigma^j(X^{{u}}(t),{u}(t))^{\top}\Phi_{xx}(X^{{u}}(t))\sigma^j(X^{{u}}(t)
,{u}(t)).
$$
We introduce the processes $(p^{\varepsilon}_{0}(\cdot),K^{\varepsilon}_{0}(\cdot))$ and $(P^{\varepsilon}_{0}(\cdot),Q^{\varepsilon}_{0}(\cdot))$ which satisfy the following asymptotic first- and second-order adjoint equations.
\begin{equation}\label{afirst}
\begin{aligned}
-dp^{\varepsilon}_{0}(t)&=\left[b_x^*(\bar{X}(t), \bar{u}(t))p^{\varepsilon}_{0}(t)+\sum_{j=1}^d \sigma_x^{j*}(\bar{X}(t), \bar{u}(t)) K_0^{{\varepsilon},j}(t)+l_x(\bar{X}(t), \bar{u}(t))\right] \mathrm{d}t-K^{\varepsilon}_{0}(t) \mathrm{d} W(t), \\
p^{\varepsilon}_{0}(\tau^{u^{\varepsilon}})&=0,\\
\end{aligned}
\end{equation}
and
\begin{equation}\label{asecond}
	\begin{aligned}
		 -dP^{\varepsilon}_{0}(t)=&\left[b_x^*(\bar{X}(t), \bar{u}(t)) P^{\varepsilon}_{0}(t)+P^{\varepsilon}_{0}(t) b_x(\bar{X}(t), \bar{u}(t))+\sum_{j=1}^d \sigma_x^{j *}(\bar{X}(t), \bar{u}(t)) P^{\varepsilon}_{0}(t) \sigma_x^j(\bar{X}(t), \bar{u}(t)) \right. \\
		& \left.+\sum_{j=1}^d \sigma_x^{j *}(\bar{X}(t), \bar{u}(t)) Q_{0}^{{\varepsilon},j}(t)+\sum_{j=1}^d Q_{0}^{{\varepsilon},j}(t) \sigma_x^j(\bar{X}(t), \bar{u}(t))+H^{0}_{x x}(\bar{X}(t), \bar{u}(t), p^{\varepsilon}_{0}(t), K^{\varepsilon}_{0}(t))\right] \mathrm{d}t\\
&-Q^{\varepsilon}_{0}(t) \mathrm{d} W(t),\\
		 P^{\varepsilon}_{0}(\tau^{{u}^{\varepsilon}})=&0,\\
	\end{aligned}
\end{equation}
which converge to the solutions of the first- and second-order adjoint equations in $L^1$,
\begin{equation}\label{afirst-1}
\begin{aligned}
-dp_{0}(t)&=\left[b_x^*(\bar{X}(t), \bar{u}(t))p_{0}(t)+\sum_{j=1}^d \sigma_x^{j*}(\bar{X}(t), \bar{u}(t)) K^j_0(t)+l_x(\bar{X}(t), \bar{u}(t))\right] \mathrm{d}t-K_{0}(t) \mathrm{d} W(t), \\
p_{0}(\tau^{\bar{u}})&=0,\\
\end{aligned}
\end{equation}
and
\begin{equation}\label{asecond-1}
	\begin{aligned}
		 -dP_{0}(t)=&\left[b_x^*(\bar{X}(t), \bar{u}(t)) P_{0}(t)+P_{0}(t) b_x(\bar{X}(t), \bar{u}(t))+\sum_{j=1}^d \sigma_x^{j *}(\bar{X}(t), \bar{u}(t)) P_{0}(t) \sigma_x^j(\bar{X}(t), \bar{u}(t)) \right. \\
		& \left.+\sum_{j=1}^d \sigma_x^{j *}(\bar{X}(t), \bar{u}(t)) Q_{0}^j(t)+\sum_{j=1}^d Q_{0}^j(t) \sigma_x^j(\bar{X}(t), \bar{u}(t))+H^{0}_{x x}(\bar{X}(t), \bar{u}(t), p_{0}(t), K_{0}(t))\right] \mathrm{d}t\\
&-Q_{0}(t) \mathrm{d} W(t),\\
		 P_{0}(\tau^{\bar{u}})=&0,\\
	\end{aligned}
\end{equation}
where
$$
H^{0}(x, v, p_{0}, k_{0})=l(x, v)+(p_{0}, b(x, v))+\sum_{j=1}^d\left(k_{0}^j, \sigma^j(x, v)\right).
$$

\begin{lemma}
\label{ale-1}
Let Assumptions \ref{ass-a} and \ref{ass-b} hold, we have
$$
\bar{h}(u,t)=\mathbb{E}[k(\tau)],
$$
where
$$
\begin{aligned}
	k(\tau)=&H
	^{0}\left(\bar{X}(\tau), u, p_{0}(\tau), K_{0}(\tau)\right)-H^{0}(\bar{X}(\tau), \bar{u}(\tau), p_{0}(\tau), K_{0}(\tau))\\
	&+\frac{1}{2}\operatorname{tr}\left[\left(\sigma\left(\bar{X}(\tau), u\right)-\sigma(\bar{X}(\tau), \bar{u}(\tau))\right)^* P_{0}(\tau)\left(\sigma\left(\bar{X}(\tau), u\right)-\sigma(\bar{X}(\tau),\bar{u}(\tau))\right)\right],
\end{aligned}
$$

\end{lemma}
\noindent {\textbf{Proof}}:
By Lemma \ref{le-1} and Assumption \ref{ass-a}, we have
\begin{equation*}
\begin{aligned}	
&\int_{0}^{\tau^{{u}^{\varepsilon}}}
\left(h^{u^{\varepsilon}}\left(t\right)-h^{\bar{u}}\left(t\right)\right)\mathrm{d}t\\
=&\int_{0}^{\tau^{{u}^{\varepsilon}}}\mathbb{E}[l\left(X^{\varepsilon}(t),
{u}^{\varepsilon}(t)\right)-l\left(\bar{X}(t),\bar{u}(t)\right)]\mathrm{d}t\\
=&\int_{0}^{\tau^{{u}^{\varepsilon}}}\mathbb{E}[l\left(\bar{X}(t)+y_{1}(t)+y_{2}(t),
{u}^{\varepsilon}(t)\right)-l\left(\bar{X}(t),\bar{u}(t)\right)]\mathrm{d}t+o(\varepsilon)\\
=&\int_{0}^{\tau^{{u}^{\varepsilon}}}\mathbb{E}[l\left(\bar{X}(t)+y_{1}(t)+y_{2}(t),
{u}^{\varepsilon}(t)\right)-l\left(\bar{X}(t),{u}^{\varepsilon}(t)\right)
+l\left(\bar{X}(t),{u}^{\varepsilon}(t)\right)-l\left(\bar{X}(t),\bar{u}(t)\right)]
\mathrm{d}t+o(\varepsilon)\\		=&\int_{0}^{\tau^{{u}^{\varepsilon}}}\mathbb{E}[l_{x}\left(\bar{X}(t)),{u}^{\varepsilon}(t)\right)
(y_{1}(t)+y_{2}(t))+\frac{1}{2}l_{xx}\left(\bar{X}(t)),{u}^{\varepsilon}(t)\right)
y_{1}(t)y_{1}(t)\\
&+l\left(\bar{X}(t),{u}^{\varepsilon}(t)\right)-l\left(\bar{X}(t),\bar{u}(t)\right)]\mathrm{d}t
+o(\varepsilon)\\		=&\int_{0}^{\tau^{{u}^{\varepsilon}}}\mathbb{E}[l_{x}\left(\bar{X}(t)),\bar{u}(t)\right)(y_{1}(t)
+y_{2}(t))+\frac{1}{2}l_{xx}\left(\bar{X}(t)),\bar{u}(t)\right)y_{1}(t)y_{1}(t)\\
&+l\left(\bar{X}(t),{u}^{\varepsilon}(t)\right)-l\left(\bar{X}(t),\bar{u}(t)\right)]
\mathrm{d}t+o(\varepsilon).
\end{aligned}
\end{equation*}
Thus
\begin{equation}
\label{aeq-3}
\begin{aligned}	
&\int_{0}^{\tau^{{u}^{\varepsilon}}}
\left(h^{u^{\varepsilon}}\left(t\right)-h^{\bar{u}}\left(t\right)\right)\mathrm{d}t\\
=&\int_{0}^{\tau^{{u}^{\varepsilon}}}\mathbb{E}[l_{x}\left(\bar{X}(t)),\bar{u}(t)\right)(y_{1}(t)
+y_{2}(t))+\frac{1}{2}l_{xx}\left(\bar{X}(t)),\bar{u}(t)\right)y_{1}(t)y_{1}(t)\\
&+l\left(\bar{X}(t),{u}^{\varepsilon}(t)\right)-l\left(\bar{X}(t),\bar{u}(t)\right)]
\mathrm{d}t+o(\varepsilon).
\end{aligned}
\end{equation}

By variational equations (\ref{eq-1}), (\ref{eq-2}) and first-order adjoint equation (\ref{afirst}), we have
\begin{equation}
	\begin{aligned}
		&\mathbb{E}\int_{0}^{\tau^{{u}^{\varepsilon}}} l_x(s) y_1(s)\mathrm{d}s \\
		=&\mathbb{E}\int _{0}^{\tau^{{u}^{\varepsilon}}} \left(p^{\varepsilon}_{0}(s), b(\bar{X}(s), u^{\varepsilon}(s))-b(\bar{X}(s), \bar{u}(s))\right) \mathrm{d}s\\
		&+\mathbb{E}\int _{0}^{\tau^{{u}^{\varepsilon}}}  \operatorname{tr}\left[K^{\varepsilon}_{0}(s)\left(\sigma\left(\bar{X}(s), u^{\varepsilon}(s)\right)-\sigma(\bar{X}(s), \bar{u}(s))\right] \mathrm{d}s\right. \\
	\end{aligned}
\end{equation}
and
\begin{equation}
	\begin{aligned}
&\mathbb{E}\int _{0}^{\tau^{{u}^{\varepsilon}}}l_x(s) y_2(s) \mathrm{d}s\\
		=&\mathbb{E}\int _{0}^{\tau^{{u}^{\varepsilon}}} \frac{1}{2}\left[\left(p^{\varepsilon}_{0}(s) b_{x x}(s)+\sum_{j=1}^d K_{0}^{\varepsilon,j}(s) \sigma_{x x}^j(s)\right) y_1(s) y_1(s)\right] \mathrm{d}s\\
		&+\mathbb{E}\int _{0}^{\tau^{{u}^{\varepsilon}}} p^{\varepsilon*}_{0}(s)\left(b_x\left(\bar{X}(s), u^{\varepsilon}(s)\right)-b_{x}(\bar{X}(s), \bar{u}(s))\right) y_{1}(s) \mathrm{d}s  \\
		&+\mathbb{E}\int _{0}^{\tau^{{u}^{\varepsilon}}} \sum_{j=1}^d K_0^{{\varepsilon},j*}(s)\left(\sigma_x^j\left(\bar{X}(s), u^{\varepsilon}(s)\right)-\sigma_x^j(\bar{X}(s), \bar{u}(s))\right) y_1(s) \mathrm{d}s .
	\end{aligned}
\end{equation}
Thus we can rewrite (\ref{aeq-3}) as
\begin{equation}
	\label{aeq-5}
	\begin{aligned}
		&\mathbb{E} \int_{0}^{\tau^{{u}^{\varepsilon}}}\left(H^{0}\left(\bar{X}(s), u^{\varepsilon}(s), p^{\varepsilon}_{0}(s), K^{\varepsilon}_{0}(s)\right)-H^{0}(\bar{X}(s), \bar{u}(s), p^{\varepsilon}_{0}(s), K^{\varepsilon}_{0}(s))\right) \mathrm{d}s\\
		&+\frac{1}{2} \mathbb{E} \int_{0}^{\tau^{{u}^{\varepsilon}}} y_1^*(s) H^{0}_{x x}(\bar{X}(s), \bar{u}(s), p^{\varepsilon}_{0}(s), K_0^{\varepsilon}(s)) y_1(s)\mathrm{d}s. \\
	\end{aligned}
\end{equation}

Again, by second-order adjoint equation (\ref{asecond}), we rewrite (\ref{aeq-3}) as
$$
\begin{aligned}
	&\mathbb{E}\int_{0}^{\tau^{{u}^{\varepsilon}}}\left(H^{0}\left(\bar{X}(s), u^{\varepsilon}(s), p^{\varepsilon}_{0}(s), K^{\varepsilon}_{0}(s)\right)-H^{0}(\bar{X}(s), \bar{u}(s), p^{\varepsilon}_{0}(s), K^{\varepsilon}_{0}(s))\right) \mathrm{d}s \\
	&+\frac{1}{2} \mathbb{E} \int_{0}^{\tau^{{u}^{\varepsilon}}} \operatorname{tr}\left[\left(\sigma\left(\bar{X}(s), u^{\varepsilon}(s)\right)-\sigma(\bar{X}(s), \bar{u}(s))\right)^* P^{\varepsilon}_{0}(s)\left(\sigma\left(\bar{X}(s), u^{\varepsilon}(s)\right)-\sigma(\bar{X}(s),\bar{u}(s))\right)\right] \mathrm{d}s\\
	=&\mathbb{E}\left[(H
	^{0}\left(\bar{X}(\tau), u, p^{\varepsilon}_{0}(\tau), K^{\varepsilon}_{0}(\tau)\right)-H^{0}(\bar{X}(\tau), \bar{u}(\tau), p^{\varepsilon}_{0}(\tau), K^{\varepsilon}_{0}(\tau)))\varepsilon\right]\\
	&+\mathbb{E}\left[\frac{1}{2}\operatorname{tr}\left[\left(\sigma\left(\bar{X}(\tau), u\right)-\sigma(\bar{X}(\tau), \bar{u}(\tau))\right)^* P^{\varepsilon}_{0}(\tau)\left(\sigma\left(\bar{X}(\tau), u\right)-\sigma(\bar{X}(\tau),\bar{u}(\tau))\right)\right]\right]\varepsilon+o(\varepsilon).
\end{aligned}
$$
Therefore,
\begin{equation}\label{avar-1}
\bar{h}(u,t):=\lim_{\varepsilon\to 0}\int_{0}^{\tau^{{u}^{\varepsilon}}}\frac{h^{u^{\varepsilon}}
\left(t\right)-h^{\bar{u}}\left(t\right)}{\varepsilon}\mathrm{d}t=\mathbb{E}\left[k(\tau)\right],\\
\end{equation}
where
$$
\begin{aligned}
	k(\tau)=&H
	^{0}\left(\bar{X}(\tau), u, p_{0}(\tau), K_{0}(\tau)\right)-H^{0}(\bar{X}(\tau), \bar{u}(\tau), p_{0}(\tau), K_{0}(\tau))\\
	&+\frac{1}{2}\operatorname{tr}\left[\left(\sigma\left(\bar{X}(\tau), u\right)-\sigma(\bar{X}(\tau), \bar{u}(\tau))\right)^* P_{0}(\tau)\left(\sigma\left(\bar{X}(\tau), u\right)-\sigma(\bar{X}(\tau),\bar{u}(\tau))\right)\right].\\
\end{aligned}
$$
\noindent This completes the proof. $\ \ \ \ \ \ \ \ \Box$
\bigskip

\begin{lemma}
\label{le-4}
Let  Assumptions \ref{ass-a} and \ref{ass-b} hold. Suppose that $h^{\bar{u}}(\tau^{\bar{u}})\neq 0$, and $h^{\bar{u}}(\cdot)$ is continuous at the point $\tau^{\bar{u}}$. We have the following results.
	
(i). If $\tau^{\bar{u}}<T$, one obtains
\begin{equation}
\begin{aligned}
&\lim_{\varepsilon\to 0}\frac{\tau^{\bar{u}}-\tau^{u^{\varepsilon}}}{\varepsilon}=
\frac{\mathbb{E}\left[k(\tau)\right]}{h^{\bar{u}}(\tau^{\bar{u}})},\\
\end{aligned}
\end{equation}
where
$$
\begin{aligned}
k(\tau)=& H
^{0}\left(\bar{X}(\tau), u, p_{0}(\tau), K_{0}(\tau)\right)-H^{0}(\bar{X}(\tau), \bar{u}(\tau), p_{0}(\tau), K_{0}(\tau))\\
&+\frac{1}{2}\operatorname{tr}\left[\left(\sigma\left(\bar{X}(\tau), u\right)-\sigma(\bar{X}(\tau), \bar{u}(\tau))\right)^* P_{0}(\tau)\left(\sigma\left(\bar{X}(\tau), u\right)-\sigma(\bar{X}(\tau),\bar{u}(\tau))\right)\right].\\
\end{aligned}
$$

	(ii). If $\inf\bigg{\{}t:\mathbb{E}[\Phi(\bar{X}(t))]\geq \alpha,\ t\in [0,T] \bigg{\}}=T$, then there exists a sequence $\varepsilon_n\to 0$ as $n\to +\infty$ such that
\begin{equation}
\begin{aligned}
&\lim_{n\to +\infty}\frac{\tau^{\bar{u}}-\tau^{u^{\varepsilon_n}}}{\varepsilon_n}
=\frac{\mathbb{E}\left[k(\tau)\right]}{h^{\bar{u}}(\tau^{\bar{u}})}\ \  \mathrm{or}\ \  0.\\
\end{aligned}
\end{equation}

(iii). If $\bigg{\{}t:\mathbb{E}[\Phi(\bar{X}(t))]\geq \alpha,\ t\in [0,T] \bigg{\}}=\varnothing$, we have
	\begin{equation}
		\lim_{\varepsilon\to 0}\frac{\tau^{\bar{u}}-\tau^{u^{\varepsilon}}}{\varepsilon}=0.
	\end{equation}	
\end{lemma}
\textbf{Proof:} We first prove case $(i)$. Notice that for $ \tau^{\bar{u}}<T$,
\begin{equation*}
	\tau^{\bar{u}}=\inf\bigg{\{}t:\mathbb{E}[\Phi({\bar{X}}(t))]\geq \alpha,\ t\in [0,T] \bigg{\}}.
\end{equation*}
By $(i)$ of Lemma \ref{le-2}, we have that $\displaystyle\lim_{\varepsilon\to 0}\left|\tau^{u^{\varepsilon}}-\tau^{\bar{u}}\right|=0$, for sufficiently small $\varepsilon>0$, it follows that
\begin{equation*}
\tau^{{u^{\varepsilon}}}=\inf\bigg{\{}t:\mathbb{E}[\Phi(X^{\varepsilon}(t))]\geq \alpha,\ t\in [0,T] \bigg{\}}.
\end{equation*}
This implies that
$
\mathbb{E}[\Phi({\bar{X}}(\tau^{\bar{u}}))]=\mathbb{E}[\Phi(X^{\varepsilon}(\tau^{u^{\varepsilon}}))]=\alpha.
$
Combining
$$
\mathbb{E}[\Phi(\bar{X}(\tau^{\bar{u}}))]=\Phi(x_0)+
\int_0^{\tau^{\bar{u}}}h^{\bar{u}}(t)\mathrm{d}t,
$$
and
$$
\mathbb{E}[\Phi(X^{\varepsilon}(\tau^{u^{\varepsilon}}))]=\Phi(x_0)
+\int_0^{\tau^{u^{\varepsilon}}}h^{u^{\varepsilon}}(t) \mathrm{d}t,
$$
we have
$$
\int_0^{\tau^{\bar{u}}}h^{\bar{u}}(t)\mathrm{d}t=\int_0^{\tau^{u^{\varepsilon}}}h^{{u^{\varepsilon}}}(t)\mathrm{d}t,
$$
and
$$
\int_{\tau^{u^{\varepsilon}}}^{\tau^{\bar{u}}}h^{\bar{u}}(t)\text{d}t
=\int_0^{\tau^{u^{\varepsilon}}}\bigg{[}h^{u^{\varepsilon}}(t)-h^{\bar{u}}(t)\bigg{]}\text{d}t.
$$
Dividing  on both sides of the above equation by $\varepsilon$, we obtain
$$
\displaystyle \lim_{\varepsilon\to 0}\frac{\int_{\tau^{u^{\varepsilon}}}^{\tau^{\bar{u}}}h^{\bar{u}}(t)\text{d}t}{\varepsilon}
=\displaystyle \lim_{\varepsilon\to 0}\int_0^{\tau^{u^{\varepsilon}}}\frac{h^{u^{\varepsilon}}(t)-h^{\bar{u}}(t)}{\varepsilon}\text{d}t.
$$

By $(i)$ of Lemma \ref{le-1}, $\displaystyle\lim_{\varepsilon\to 0}\left|\tau^{u^{\varepsilon}}-\tau^{\bar{u}}\right|=0$, we have
\begin{equation}
\label{le2-e1}
\begin{array}
		[c]{rl}%
		& \displaystyle \lim_{\varepsilon\to 0}\frac{\int_{\tau^{u^{\varepsilon}}}^{\tau^{\bar{u}}}h^{\bar{u}}(t)\text{d}t}{\varepsilon}
		= \displaystyle \lim_{\varepsilon\to 0}\frac{\tau^{u^{\varepsilon}}-\tau^{\bar{u}}}{\varepsilon}\bigg{[}h^{\bar{u}}(\tau^{\bar{u}})
		+o(1)\bigg{]},\\
\end{array}
\end{equation}
where $o(1)$ converges to $0$ as $\left|\tau^{u^{\varepsilon}}-\tau^{\bar{u}}\right|\to 0$. By formula $(\ref{avar-1})$, it follows that
\begin{equation}\label{le2-ee1}
\begin{aligned}
&\lim_{\varepsilon\to 0}\int_{0}^{\tau^{{u}^{\varepsilon}}}\frac{h^{u^{\varepsilon}}\left(t\right)-h^{\bar{u}}\left(t\right)}{\varepsilon}dt\\
=&\mathbb{E}\bigg[H^{0}\left(\bar{X}(\tau), u, p_{0}(\tau), K_{0}(\tau)\right)-H^{0}(\bar{X}(\tau), \bar{u}(\tau), p_{0}(\tau), K_{0}(\tau))\\
&+\frac{1}{2}\operatorname{tr}\left[\left(\sigma\left(\bar{X}(\tau), u\right)-\sigma(\bar{X}(\tau), \bar{u}(\tau))\right)^* P_{0}(\tau)\left(\sigma\left(\bar{X}(\tau), u\right)-\sigma(\bar{X}(\tau),\bar{u}(\tau))\right)\right]\bigg].\\
\end{aligned}
\end{equation}

Combining equations (\ref{le2-e1}) and (\ref{le2-ee1}), we obtain
\begin{equation}
\begin{aligned}
		&\lim_{\varepsilon\to 0}\frac{\tau^{\bar{u}}-\tau^{u^{\varepsilon}}}{\varepsilon} \\
=&\frac{1}{h^{\bar{u}}(\tau^{\bar{u}})}\mathbb{E}\left[H
		^{0}\left(\bar{X}(\tau), u, p_{0}(\tau), K_{0}(\tau)\right)-H^{0}(\bar{X}(\tau), \bar{u}(\tau), p_{0}(\tau), K_{0}(\tau))\right.\\
		&\left.+\frac{1}{2}\operatorname{tr}\left[\left(\sigma\left(\bar{X}(\tau), u\right)-\sigma(\bar{X}(\tau), \bar{u}(\tau))\right)^* P_{0}(\tau)\left(\sigma\left(\bar{X}(\tau), u\right)-\sigma(\bar{X}(\tau),\bar{u}(\tau))\right)\right]\right].\\
\end{aligned}
\end{equation}

{ Secondly, we consider the case $(ii)$. Notice that $\inf\bigg{\{}t:\mathbb{E}[\Phi(\bar{X}(t))]\geq \alpha,\ t\in [0,T] \bigg{\}}=T$, if there exists sequence $\varepsilon_n\to 0$ as $n\to +\infty$ such that
	$$
	\tau^{{u^{\varepsilon_n}}}=\inf\bigg{\{}t:\mathbb{E}[\Phi(X^{\varepsilon_n}(t))]\geq \alpha,\ t\in [0,T] \bigg{\}}<T.
	$$
	Similar with the proof of case $(i)$, by the case $(ii)$ of Lemma \ref{le-2}, we can obtain
	\begin{equation*}
	\begin{aligned}
			&\lim_{n\to +\infty}\frac{\tau^{\bar{u}}-\tau^{u^{\varepsilon_n}}}{\varepsilon_n} \\
=&\frac{1}{h^{\bar{u}}(\tau^{\bar{u}})}\mathbb{E}\bigg[H
			^{0}\left(\bar{X}(\tau), u, p_{0}(\tau), K_{0}(\tau)\right)-H^{0}(\bar{X}(\tau), \bar{u}(\tau), p_{0}(\tau), K_{0}(\tau))\\
			&+\frac{1}{2}\operatorname{tr}\left[\left(\sigma\left(\bar{X}(\tau), u\right)-\sigma(\bar{X}(\tau), \bar{u}(\tau))\right)^* P_{0}(\tau)\left(\sigma\left(\bar{X}(\tau), u\right)-\sigma(\bar{X}(\tau),\bar{u}(\tau))\right)\right]\bigg].\\
	\end{aligned}
	\end{equation*}
	If there exists a sequence $\varepsilon_n\to 0$ as $n\to +\infty$ such that
	$$
	\inf\bigg{\{}t:\mathbb{E}[\Phi(X^{\varepsilon_n}(t))]\geq \alpha,\ t\in [0,T] \bigg{\}}=+\infty,
	$$
	then, $\tau^{{u^{\varepsilon_n}}}=T$, and
	\begin{equation*}
		\lim_{n\to +\infty}\frac{\tau^{\bar{u}}-\tau^{u^{\varepsilon_n}}}{\varepsilon_n}=0.
	\end{equation*}
	Thus, the case $(ii)$ is right.
	
	In the end, we consider case $(iii)$. Notice that $\bigg{\{}t:\mathbb{E}[\Phi(\bar{X}(t))]\geq \alpha,\ t\in [0,T] \bigg{\}}=\varnothing$, thus for sufficiently small $\varepsilon$,
	$$
	\inf\bigg{\{}t:\mathbb{E}[\Phi(X^{\varepsilon}(t))]\geq \alpha,\ t\in [0,T] \bigg{\}}=+\infty,
	$$
	and $\tau^{{u^{\varepsilon}}}=T$, which shows that case $(iii)$ is right.
}
\noindent This completes the proof. $\ \ \ \ \ \ \ \ \Box$
\bigskip

Next, we derive the variational inequality from the fact that
$$J\left(u^{\varepsilon}\left(\cdot\right)\right)-J\left(\bar{u}\left(\cdot\right)\right)\geq 0.$$
Applying Lemmas \ref{le-3} and \ref{le-4}, we can obtain the following result.
  \begin{lemma}
  	\label{le-5}
  Let Assumptions \ref{ass-a} and \ref{ass-b} hold, and suppose that $h^{\bar{u}}(\tau^{\bar{u}})\neq 0$ and $h^{\bar{u}}(\cdot)$ is continuous at the point $\tau^{\bar{u}}$. We have the following results.

   (i). If $\tau^{\bar{u}}<T$, one obtains
  \begin{equation}\label{vi-4}
  \begin{aligned}  				&\mathbb{E}\int_{0}^{\tau^{\bar{u}}}\big{[}f_{x}(\bar{X}(t),\bar{u}(t))(y_{1}(t)+y_{2}(t))
  +\frac{1}{2}f_{xx}(\bar{X}(t),\bar{u}(t))y_{1}(t)y_{1}(t)\big{]}\mathrm{d}t\\
  &+\mathbb{E}\int_{0}^{\tau^{\bar{u}}}\big{[}f(\bar{X}(t),u^{\varepsilon}(t))-f(\bar{X}(t),\bar{u}(t))
  \big{]}\mathrm{d}t			+\mathbb{E}\left[\varepsilon\frac{k(\tau)}{h^{\bar{u}}(\tau^{\bar{u}})}R(\tau^{\bar{u}})\right]\\
  &+\mathbb{E}\left[g_{x}(\bar{X}(\tau^{\bar{u}}))(y_{1}(\tau^{\bar{u}})+y_{2}(\tau^{\bar{u}}))\right]+\frac{1}{2}\mathbb{E}\left[g_{xx}(\bar{X}(\tau^{\bar{u}}))y_{1}(\tau^{\bar{u}})y_{1}(\tau^{\bar{u}})\right]\geq o (\varepsilon),\\
  \end{aligned}
  \end{equation}
where
  $$
  \begin{aligned}
  			R(\tau^{\bar{u}})=&f\left ( \bar{X}\left ( \tau^{\bar{u}} \right )+y_{1}\left ( \tau^{\bar{u}} \right ) +y_{2}\left ( \tau^{\bar{u}} \right ),\bar{u}\left ( \tau^{\bar{u}} \right )\right )+g_{x}(\bar{X}(\tau ^{\bar{u}} )+y_{1}( \tau ^{\bar{u}}) +y_{2}(\tau ^{\bar{u}}))B(\bar{X}(\tau^{\bar{u}}),\bar{u}(\tau ^{\bar{u}}))\\
  			&+\frac{1}{2}g_{xx}\left(\bar{X}(\tau ^{\bar{u}})+y_{1}(\tau ^{\bar{u}})+y_{2}(\tau ^{\bar{u}})\right)A(\bar{X}(\tau^{\bar{u}}),\bar{u}(\tau ^{\bar{u}})),\\
  			k(\tau)=& H
  			^{0}\left(\bar{X}(\tau), u, p_{0}(\tau), K_{0}(\tau)\right)-H^{0}(\bar{X}(\tau), \bar{u}(\tau), p_{0}(\tau), K_{0}(\tau))\\
  			&+\frac{1}{2}\operatorname{tr}\left[\left(\sigma\left(\bar{X}(\tau), u\right)-\sigma(\bar{X}(\tau), \bar{u}(\tau))\right)^* P_{0}(\tau)\left(\sigma\left(\bar{X}(\tau), u\right)-\sigma(\bar{X}(\tau),\bar{u}(\tau))\right)\right].\\
  \end{aligned}
  $$	

  (ii). If $\inf\bigg{\{}t:\mathbb{E}[\Phi(\bar{X}(t))]\geq \alpha,\ t\in [0,T] \bigg{\}}=T$, one obtains
  \begin{equation}
  \begin{aligned}  			&\mathbb{E}\int_{0}^{\tau^{\bar{u}}}\big{[}f_{x}(\bar{X}(t),\bar{u}(t))(y_{1}(t)+y_{2}(t))
  +\frac{1}{2}f_{xx}(\bar{X}(t),\bar{u}(t))y_{1}(t)y_{1}(t)\big{]}\mathrm{d}t\\
  			&+\mathbb{E}\int_{0}^{\tau^{\bar{u}}}\big{[}f(\bar{X}(t),u^{\varepsilon}(t))-f(\bar{X}(t),\bar{u}(t))
  \big{]}\mathrm{d}t 			+\mathbb{E}\left[\varepsilon\frac{k(\tau)}{h^{\bar{u}}(\tau^{\bar{u}})}R(\tau^{\bar{u}})\right]\\
  			&+\mathbb{E}\left[g_{x}(\bar{X}(\tau^{\bar{u}}))(y_{1}(\tau^{\bar{u}})+y_{2}(\tau^{\bar{u}}))\right]+\frac{1}{2}\mathbb{E}\left[g_{xx}(\bar{X}(\tau^{\bar{u}}))y_{1}(\tau^{\bar{u}})y_{1}(\tau^{\bar{u}})\right]\geq o (\varepsilon),\\
  \end{aligned}
  \end{equation}
  	    or
  \begin{equation}
  \begin{aligned}  				&\mathbb{E}\int_{0}^{\tau^{\bar{u}}}\big{[}f_{x}\left(\bar{X}\left(t\right),\bar{u}\left(t\right)\right)\left(y_{1}\left(t\right)+y_{2}\left(t\right)\right)
  +\frac{1}{2}f_{xx}\left(\bar{X}\left(t\right),\bar{u}\left(t\right)\right)y_{1}\left(t\right)y_{1}\left(t\right)\big{]}\mathrm{d}t\\
  &+\mathbb{E} \int_{0}^{\tau^{\bar{u}}}\big{[}f\left(\bar{X}(t), u^\varepsilon(t)\right)-f(\bar{X}(t), \bar{u}(t))\big{]}\mathrm{d}t  				+\mathbb{E}\left[g_x\left(\bar{X}\left(\tau^{\bar{u}}\right)\right)\left(y_1
  \left(\tau^{\bar{u}}\right)+y_2\left(\tau^{\bar{u}}\right)\right)\right] \\
  &+\frac{1}{2} \mathbb{E}\left[g_{x x}\left(\bar{X}\left(\tau^{\bar{u}}\right)\right) y_1\left(\tau^{\bar{u}}\right) y_1\left(\tau^{\bar{u}}\right)\right] \geqslant  o (\varepsilon).
  \end{aligned}
  \end{equation}
  		
  (iii). If $\bigg{\{}t:\mathbb{E}[\Phi(\bar{X}(t))]\geq \alpha,\ t\in [0,T] \bigg{\}}=\varnothing$, we have
  \begin{equation}
  \begin{aligned}
  &\mathbb{E}\int_{0}^{\tau^{\bar{u}}}\big{[}f_{x}\left(\bar{X}\left(t\right),\bar{u}\left(t\right)\right)\left(y_{1}\left(t\right)+y_{2}\left(t\right)\right)
  +\frac{1}{2}f_{xx}\left(\bar{X}\left(t\right),\bar{u}\left(t\right)\right)y_{1}\left(t\right)y_{1}\left(t\right)\big{]}\mathrm{d}t\\
  &+\mathbb{E} \int_{0}^{\tau^{\bar{u}}}\big{[}f\left(\bar{X}(t), u^\varepsilon(t)\right)-f(\bar{X}(t), \bar{u}(t))\big{]}\mathrm{d}t  				+\mathbb{E}\left[g_x\left(\bar{X}\left(\tau^{\bar{u}}\right)\right)\left(y_1
  \left(\tau^{\bar{u}}\right)+y_2\left(\tau^{\bar{u}}\right)\right)\right] \\
  &+\frac{1}{2} \mathbb{E}\left[g_{x x}\left(\bar{X}\left(\tau^{\bar{u}}\right)\right) y_1\left(\tau^{\bar{u}}\right) y_1\left(\tau^{\bar{u}}\right)\right] \geqslant  o (\varepsilon).
  			\end{aligned}
  		\end{equation}
  	\end{lemma}
\textbf{Proof:} We first prove case $(i)$. Since $\left(\bar{X}\left(\cdot\right),\bar{u}\left(\cdot\right)\right)$ is an optimal state of the control problem, we have $J\left(u^{\varepsilon}\left(\cdot\right)\right)-J\left(\bar{u}\left(\cdot\right)\right)\geqslant 0$, i.e.
  	$$\mathbb{E}\left[\int_{0}^{\tau^{u^{\varepsilon}}}f\left(X^{\varepsilon}(t)
  ,u^{\varepsilon}(t)\right)\mathrm{d}t-g(X^{\varepsilon}
  (\tau^{u^\varepsilon}))\right]-\mathbb{E}
  \left[\int_{0}^{\tau^{\bar{u}}}f\left(\bar{X}(t),\bar{u}
  (t)\right)\mathrm{d}t-g(\bar{X}(\tau^{\bar{u}}))\right] \geqslant 0.$$
Thus, from Lemma \ref{le-3}, Lemma \ref{le-4} and the following  moment inequality
  	$$\underset{0\leqslant t \leqslant T}\sup{\mathbb{E}\left | y_{1}(t) \right | ^2}\leqslant C\varepsilon,$$
  	$$\underset{0\leqslant t \leqslant T}\sup{\mathbb{E}\left | y_{2}(t) \right | ^2}\leqslant C\varepsilon^{2},$$
  	$$\underset{0\leqslant t \leqslant T}\sup{\mathbb{E}\left | y_{1}(t) \right | ^4}\leqslant C\varepsilon^{2},$$
  	$$\underset{0\leqslant t \leqslant T}\sup{\mathbb{E}\left | y_{2}(t) \right | ^4}\leqslant C\varepsilon^{4},$$
  	we have
  	\begin{equation}
  		\begin{aligned}
  			0 \leqslant &\mathbb{E} \int_{0}^{\tau^{u^\varepsilon }}\big{[}f\left ( \bar{X}\left ( t \right )+y_{1}\left ( t \right ) +y_{2}\left ( t \right ),u^{\varepsilon }\left ( t \right ) \right )-f\left (\bar{X}\left ( t \right ),\bar{u}\left ( t \right )   \right )\big{]}\mathrm{d}t\\
  			&+\mathbb{E}\left [ g\left ( \bar{X}\left ( \tau^{u^\varepsilon } \right )+y_{1}\left ( \tau^{u^\varepsilon } \right ) +y_{2}\left ( \tau^{u^\varepsilon }\right ) \right ) -g\left ( \bar{X} \left ( \tau ^{\bar{u}} \right )  \right )  \right ]+o (\varepsilon) \\
  			=&\mathbb{E} \int_{0}^{\tau^{\bar{u}}}\big{[}f\left ( \bar{X}\left ( t \right )+y_{1}\left ( t \right ) +y_{2}\left ( t \right ),u^{\varepsilon }\left ( t \right ) \right )-f\left (\bar{X}\left ( t \right )+y_{1}\left ( t \right ) +y_{2}\left ( t \right ),\bar{u}\left ( t \right )   \right )\big{]}\mathrm{d}t\\
  			&+\mathbb{E} \int_{0}^{\tau^{\bar{u}}}\big{[}f\left (\bar{X}\left ( t \right )+y_{1}\left ( t \right ) +y_{2}\left ( t \right ),\bar{u}\left ( t \right )   \right )-f\left (\bar{X}\left ( t \right ),\bar{u}\left ( t \right )\right )\big{]}\mathrm{d}t\\
  			&+\mathbb{E} \int_{\tau^{\bar{u}}}^{\tau^{u^\varepsilon }}f\left ( \bar{X}\left ( t \right )+y_{1}\left ( t \right ) +y_{2}\left ( t \right ),u^{\varepsilon }\left ( t \right ) \right )\mathrm{d}t\\
  			&+\mathbb{E}\left[ g\left ( \bar{X}\left ( \tau^{u^\varepsilon } \right )+y_{1}\left ( \tau^{u^\varepsilon } \right ) +y_{2}\left ( \tau^{u^\varepsilon }\right ) \right ) - g\left ( \bar{X}\left (\tau ^{\bar{u}} \right )+y_{1}\left ( \tau ^{\bar{u}} \right ) +y_{2}\left (\tau ^{\bar{u}}\right ) \right )\right]\\
  			&+\mathbb{E}\left[g\left ( \bar{X}\left (\tau ^{\bar{u}} \right )+y_{1}\left ( \tau ^{\bar{u}} \right ) +y_{2}\left (\tau ^{\bar{u}}\right ) \right )-g\left ( \bar{X} \left ( \tau ^{\bar{u}} \right )  \right )\right]+o (\varepsilon).
  		\end{aligned}
    \end{equation}

 Denote by
 \begin{equation*}
 \begin{aligned}
 &I_{1}:=\varepsilon^{-1}\mathbb{E} \int_{0}^{\tau^{\bar{u}}}\big{[}f\left ( \bar{X}\left ( t \right )+y_{1}\left ( t \right ) +y_{2}\left ( t \right ),u^{\varepsilon }\left ( t \right ) \right )-f\left (\bar{X}\left ( t \right )+y_{1}\left ( t \right ) +y_{2}\left ( t \right ),\bar{u}\left ( t \right )   \right )\big{]}\mathrm{d}t,\\
 &I_{2}:=\varepsilon^{-1}\mathbb{E} \int_{0}^{\tau^{\bar{u}}}\big{[}f\left (\bar{X}\left ( t \right )+y_{1}\left ( t \right ) +y_{2}\left ( t \right ),\bar{u}\left ( t \right )   \right )-f\left (\bar{X}\left ( t \right ),\bar{u}\left ( t \right )\right )\big{]}\mathrm{d}t,\\
 &I_{3}:=\varepsilon^{-1}\mathbb{E} \int_{\tau^{\bar{u}}}^{\tau^{u^\varepsilon }}f\left ( \bar{X}\left ( t \right )+y_{1}\left ( t \right ) +y_{2}\left ( t \right ),u^{\varepsilon }\left ( t \right ) \right )\mathrm{d}t,\\
 &I_{4}:=\varepsilon^{-1}\mathbb{E}\left[ g\left ( \bar{X}\left ( \tau^{u^\varepsilon } \right )+y_{1}\left ( \tau^{u^\varepsilon } \right ) +y_{2}\left ( \tau^{u^\varepsilon }\right ) \right ) - g\left ( \bar{X}\left (\tau ^{\bar{u}} \right )+y_{1}\left ( \tau ^{\bar{u}} \right ) +y_{2}\left (\tau ^{\bar{u}}\right ) \right )\right],\\
& I_{5}:=\varepsilon^{-1}\mathbb{E}\left[g\left ( \bar{X}\left (\tau ^{\bar{u}} \right )+y_{1}\left ( \tau ^{\bar{u}} \right ) +y_{2}\left (\tau ^{\bar{u}}\right ) \right )-g\left ( \bar{X} \left ( \tau ^{\bar{u}} \right )  \right )\right].
\end{aligned}
\end{equation*}
 	
 We first consider
 	\begin{equation}
 		\begin{aligned}
 	I_{1}:=&\varepsilon^{-1}\mathbb{E}\int_{0}^{\tau^{\bar{u}}}\big{[}f\left ( \bar{X}\left ( t \right ),u^{\varepsilon }\left ( t \right ) \right )-f\left ( \bar{X}\left ( t \right ),\bar{u}\left ( t \right ) \right )\\
 	&+\left(f_{x}\left ( \bar{X}\left ( t \right ),u^{\varepsilon }\left ( t \right ) \right )\right)-f_{x}\left ( \bar{X}\left ( t \right ),\bar{u}\left ( t \right ) \right )\left(y_{1}\left ( t \right ) +y_{2}\left ( t \right )\right)\\
 	&+\frac{1}{2}\left(f_{xx}\left ( \bar{X}\left ( t \right ),u^{\varepsilon }\left ( t \right ) \right )\right)-f_{xx}\left ( \bar{X}\left ( t \right ),\bar{u}\left ( t \right ) \right )y_{1}\left ( t \right )y_{1}\left ( t \right )\big{]}\mathrm{d}t+o (1)\\
 	=&\varepsilon^{-1}\mathbb{E}\int_{0}^{\tau^{\bar{u}}}\big{[}f\left ( \bar{X}\left ( t \right ),u^{\varepsilon }\left ( t \right ) \right )-f\left ( \bar{X}\left ( t \right ),\bar{u}\left ( t \right ) \right )\big{]}\mathrm{d}t+o(1).\\
 		\end{aligned} 	
 	\end{equation}

Similarly, we obtain
\begin{equation}
	\begin{aligned}
		I_{2}:&=\varepsilon^{-1}\mathbb{E}\int_{0}^{\tau^{\bar{u}}}\big{[}f_{x}\left ( \bar{X}\left ( t \right ),\bar{u}\left ( t \right )\right ) \left(y_{1}\left ( t \right ) +y_{2}\left ( t \right )\right)+\frac{1}{2}f_{xx}\left ( \bar{X}\left ( t \right ),\bar{u}\left ( t \right ) \right )y_{1}\left ( t \right )y_{1}\left ( t \right )\big{]}\mathrm{d}t+o(1).\\
	\end{aligned} 	
\end{equation}

Applying Lemma \ref{le-2} and  Assumption \ref{ass-a}, we have
\begin{equation}
	\begin{aligned}
	I_{3}:&=\varepsilon^{-1}\mathbb{E} \int_{\tau^{\bar{u}}}^{\tau^{u^\varepsilon }}f\left ( \bar{X}\left ( t \right )+y_{1}\left ( t \right ) +y_{2}\left ( t \right ),u^{\varepsilon }\left ( t \right ) \right )\mathrm{d}t\\
	&=\mathbb{E} \left[\frac{k(\tau)}{h^{\bar{u}}(\tau^{\bar{u}})}f\left ( \bar{X}\left ( \tau^{\bar{u}} \right )+y_{1}\left ( \tau^{\bar{u}} \right ) +y_{2}\left ( \tau^{\bar{u}} \right ),\bar{u}\left ( \tau^{\bar{u}} \right )\right )\right]+o (1).\\
	\end{aligned}	
\end{equation}

Consider the term $I_{4}$. Applying It\^{o}'s formula to $g(\cdot)$, we have
\begin{equation}
	\begin{aligned}
		I_{4}=&\varepsilon^{-1}\mathbb{E}\left[ g\left ( \bar{X} (\tau^{u^\varepsilon })+y_{1} ( \tau^{u^\varepsilon }) +y_{2}( \tau^{u^\varepsilon }) \right ) - g\left ( \bar{X}(\tau ^{\bar{u}} )+y_{1}( \tau ^{\bar{u}}) +y_{2}(\tau ^{\bar{u}}) \right )\right]\\
		=&\varepsilon^{-1}\mathbb{E}\int_{\tau^{\bar{u}}}^{\tau^{u^{\varepsilon}}}g_{x}\left(\bar{X}(t)+y_{1}(t)+y_{2}(t)\right)\mathrm{d}\left(\bar{X}(t)+y_{1}(t)+y_{2}(t)\right)\\
		&+\varepsilon^{-1}\mathbb{E}\int_{\tau^{\bar{u}}}^{\tau^{u^{\varepsilon}}}\frac{1}{2}g_{xx}\left(\bar{X}(t)+y_{1}(t)+y_{2}(t)\right)\mathrm{d}\left(\left \langle \bar{X} \right \rangle _{t}+\left \langle y_{1} \right \rangle _{t}+\left \langle y_{2} \right \rangle _{t}\right)\\
=&\mathbb{E}\bigg\{\frac{k(\tau)}{h^{\bar{u}}(\tau^{\bar{u}})}[g_{x}(\bar{X}(\tau ^{\bar{u}} )+y_{1}( \tau ^{\bar{u}}) +y_{2}(\tau ^{\bar{u}}))B(\bar{X}(\tau^{\bar{u}}),\bar{u}(\tau ^{\bar{u}}))\\
&+\frac{1}{2}g_{xx}\left(\bar{X}(\tau ^{\bar{u}})+y_{1}(\tau ^{\bar{u}})+y_{2}(\tau ^{\bar{u}})\right)A(\bar{X}(\tau^{\bar{u}}),\bar{u}(\tau ^{\bar{u}}))]\bigg\}+o (1),
	\end{aligned}	
\end{equation}
where
$$
\begin{aligned}
B(\bar{X}(\tau^{\bar{u}}),\bar{u}(\tau ^{\bar{u}}))=&b(\bar{X}(\tau^{\bar{u}}),\bar{u}(\tau ^{\bar{u}}))+b_{x}(\bar{X}(\tau^{\bar{u}}),\bar{u}(\tau ^{\bar{u}}))(y_{1}(\tau^{\bar{u}})+y_{2}(\tau^{\bar{u}}))\\
&+\frac{1}{2}b_{xx}(\bar{X}(\tau^{\bar{u}}),\bar{u}(\tau ^{\bar{u}}))y_{1}(\tau^{\bar{u}})y_{1}(\tau^{\bar{u}}),\\
A(\bar{X}(\tau^{\bar{u}}),\bar{u}(\tau ^{\bar{u}}))=&\sigma^{2}(\bar{X}(\tau^{\bar{u}}),\bar{u}(\tau ^{\bar{u}}))+\sigma^{2}_{x}(\bar{X}(\tau^{\bar{u}}),\bar{u}(\tau ^{\bar{u}}))y_{1}(\tau^{\bar{u}})y_{1}(\tau^{\bar{u}}).\\
\end{aligned}
$$

For $I_{5}$, we have
\begin{equation}
	I_{5}=\varepsilon^{-1}\mathbb{E}[g_{x}(\bar{X}(\tau ^{\bar{u}}))(y_{1}( \tau ^{\bar{u}}) +y_{2}(\tau ^{\bar{u}}))+\frac{1}{2}g_{xx}(\bar{X}(\tau^{\bar{u}}))y_{1}(\tau^{\bar{u}})y_{1}(\tau^{\bar{u}})]
+o(1).
\end{equation}
Therefore, we have
\begin{equation}
\begin{aligned}	&\mathbb{E}\int_{0}^{\tau^{\bar{u}}}\big{[}f_{x}(\bar{X}(t),\bar{u}(t))(y_{1}(t)+y_{2}(t))
+\frac{1}{2}f_{xx}(\bar{X}(t),\bar{u}(t))y_{1}(t)y_{1}(t)\big{]}\mathrm{d}t\\
	&+\mathbb{E}\int_{0}^{\tau^{\bar{u}}}\big{[}f(\bar{X}(t),u^{\varepsilon}(t))-f(\bar{X}(t),\bar{u}(t))\big{]}\mathrm{d}t +\mathbb{E}\left[\varepsilon\frac{k(\tau)}{h^{\bar{u}}(\tau^{\bar{u}})}R(\tau^{\bar{u}})\right]\\
	&+\mathbb{E}\left[g_{x}(\bar{X}(\tau^{\bar{u}}))(y_{1}(\tau^{\bar{u}})+y_{2}(\tau^{\bar{u}}))\right]+\frac{1}{2}\mathbb{E}\left[g_{xx}(\bar{X}(\tau^{\bar{u}}))y_{1}(\tau^{\bar{u}})y_{1}(\tau^{\bar{u}})\right]\geq o (\varepsilon),\\
\end{aligned}
\end{equation}
where
$$
\begin{aligned}
R(\tau^{\bar{u}})=&f\left ( \bar{X}\left ( \tau^{\bar{u}} \right )+y_{1}\left ( \tau^{\bar{u}} \right ) +y_{2}\left ( \tau^{\bar{u}} \right ),\bar{u}\left ( \tau^{\bar{u}} \right )\right )+g_{x}(\bar{X}(\tau ^{\bar{u}} )+y_{1}( \tau ^{\bar{u}}) +y_{2}(\tau ^{\bar{u}}))B(\bar{X}(\tau^{\bar{u}}),\bar{u}(\tau ^{\bar{u}}))\\
&+\frac{1}{2}g_{xx}\left(\bar{X}(\tau ^{\bar{u}})+y_{1}(\tau ^{\bar{u}})+y_{2}(\tau ^{\bar{u}})\right)A(\bar{X}(\tau^{\bar{u}}),\bar{u}(\tau ^{\bar{u}})).\\
\end{aligned}
$$

Similar with the proof of case (i), by case (ii) of Lemma \ref{le-2}, one obtains
\begin{equation}
	\begin{aligned}
		&\mathbb{E}\int_{0}^{\tau^{\bar{u}}}\big{[}f_{x}(\bar{X}(t),\bar{u}(t))(y_{1}(t)+y_{2}(t))
+\frac{1}{2}f_{xx}(\bar{X}(t),\bar{u}(t))y_{1}(t)y_{1}(t)\big{]}\mathrm{d}t\\
		&+\mathbb{E}\int_{0}^{\tau^{\bar{u}}}\big{[}f(\bar{X}(t),u^{\varepsilon}(t))-f(\bar{X}(t),\bar{u}(t))\big{]}
\mathrm{d}t		+\mathbb{E}\left[\varepsilon\frac{k(\tau)}{h^{\bar{u}}(\tau^{\bar{u}})}R(\tau^{\bar{u}})\right]\\
		&+\mathbb{E}\left[g_{x}(\bar{X}(\tau^{\bar{u}}))(y_{1}(\tau^{\bar{u}})+y_{2}(\tau^{\bar{u}}))\right]+\frac{1}{2}\mathbb{E}\left[g_{xx}(\bar{X}(\tau^{\bar{u}}))y_{1}(\tau^{\bar{u}})y_{1}(\tau^{\bar{u}})\right]\geq o (\varepsilon).\\
	\end{aligned}
\end{equation}

If $$\inf\bigg{\{}t:\mathbb{E}[\Phi(X^{u^{\varepsilon}}(t))]\geq \alpha,\ t\in [0,T] \bigg{\}}=+\infty,$$
then $\tau^{u^{\varepsilon}}=T$, we can obtain
\begin{equation}
	\begin{aligned}
		&\mathbb{E}\int_{0}^{\tau^{\bar{u}}}\big{[}f_{x}\left(\bar{X}\left(t\right),\bar{u}\left(t\right)\right)\left(y_{1}\left(t\right)+y_{2}\left(t\right)\right)
+\frac{1}{2}f_{xx}\left(\bar{X}\left(t\right),\bar{u}\left(t\right)\right)y_{1}\left(t\right)y_{1}\left(t\right)\big{]}\mathrm{d}t\\
		&+\mathbb{E} \int_{0}^{\tau^{\bar{u}}}\big{[}f_{x}\left(\bar{X}(t), u^\varepsilon(t)\right)-f(\bar{X}(t), \bar{u}(t))\big{]}\mathrm{d}t \\
		&+\mathbb{E}\left[g_x\left(\bar{X}\left(\tau^{\bar{u}}\right)\right)\left(y_1\left(\tau^{\bar{u}}\right)+y_2\left(\tau^{\bar{u}}\right)\right)\right]+\frac{1}{2} \mathbb{E}\left[g_{x x}\left(\bar{X}\left(\tau^{\bar{u}}\right)\right) y_1\left(\tau^{\bar{u}}\right) y_1\left(\tau^{\bar{u}}\right)\right] \geqslant  o (\varepsilon).
	\end{aligned}
\end{equation}

In the end, for the case (iii), by the case (iii) of Lemma \ref{le-2}, for sufficiently small $\varepsilon$, we have
$$
\inf\bigg{\{}t:\mathbb{E}[\Phi(X^{u^{\varepsilon}}(t))]\geq \alpha,\ t\in [0,T] \bigg{\}}=+\infty, $$
thus, $\tau^{\bar{u}}=\tau^{u^{\varepsilon}}=T$, we have
\begin{equation}
	\begin{aligned}
		&\mathbb{E}\int_{0}^{\tau^{\bar{u}}}\big{[}f_{x}\left(\bar{X}\left(t\right),\bar{u}\left(t\right)\right)\left(y_{1}\left(t\right)+y_{2}\left(t\right)\right)
+\frac{1}{2}f_{xx}\left(\bar{X}\left(t\right),\bar{u}\left(t\right)\right)y_{1}\left(t\right)y_{1}\left(t\right)\big{]}\mathrm{d}t\\
		&+\mathbb{E} \int_{0}^{\tau^{\bar{u}}}\big{[}f_{x}\left(\bar{X}(t), u^\varepsilon(t)\right)-f(\bar{X}(t), \bar{u}(t))\big{]}\mathrm{d}t \\
		&+\mathbb{E}\left[g_x\left(\bar{X}\left(\tau^{\bar{u}}\right)\right)\left(y_1\left(\tau^{\bar{u}}\right)+y_2\left(\tau^{\bar{u}}\right)\right)\right]+\frac{1}{2} \mathbb{E}\left[g_{x x}\left(\bar{X}\left(\tau^{\bar{u}}\right)\right) y_1\left(\tau^{\bar{u}}\right) y_1\left(\tau^{\bar{u}}\right)\right] \geqslant  o (\varepsilon).
	\end{aligned}
\end{equation}
\noindent This completes the proof. $\ \ \ \ \ \ \ \ \Box$

\section{Stochastic Maximum Principle}

To establish the stochastic maximum principle under a varying terminal time for non-convex control domain, we introduce the first- and second-order adjoint processes for variational equations (\ref{eq-1}) and (\ref{eq-2}). With these processes, we can derive the variational inequality by Lemma \ref{le-5}. The following analysis is analogous to Peng \cite{P90}. For the convenience of readers and the completeness of the paper, we provide the main steps of the proof under varying terminal time.

For the case (i) of Lemma \ref{le-5}, the linear terms  in the inequality (\ref{vi-4}) can be treated in the following way. For simplicity, we let
$r_{x}(t)=r_{x}(\bar{X}(t),\bar{u}(t)),r_{xx}(t)=r_{xx}(\bar{X}(t),\bar{u}(t))$, for $r=b,\sigma, f, g$. Note that, $\tau^{\bar{u}}$ is a constant. Based on the first- and second-order adjoint equations introduced in \cite{P90}, we introduce the processes $(p(\cdot),K(\cdot))$ and $(P(\cdot),Q(\cdot))$ those are the solutions of the following first- and second-order adjoint equations,
 \begin{equation}\label{first-or}
 \begin{aligned}
 	-\mathrm{d}p(t)&=\left[b_x^*(\bar{X}(t), \bar{u}(t))p(t)+\sum_{j=1}^d \sigma_x^{j*}(\bar{X}(t), \bar{u}(t)) K_j(t)+f_x(\bar{X}(t), \bar{u}(t))\right] \mathrm{d}t-K(t) \mathrm{d} W(t), \\
 	p(\tau^{\bar{u}})&=g_x(\bar{X}(\tau^{\bar{u}})),\\
 \end{aligned}
 \end{equation}
 and
\begin{equation}\label{second-or}
	\begin{aligned}
		 -\mathrm{d}P(t)=&\left[b_x^*(\bar{X}(t), \bar{u}(t)) P(t)+P(t) b_x(\bar{X}(t), \bar{u}(t))+\sum_{j=1}^d \sigma_x^{j *}(\bar{X}(t), \bar{u}(t)) P(t) \sigma_x^j(\bar{X}(t), \bar{u}(t)) \right. \\
		& \left.+\sum_{j=1}^d \sigma_x^{j *}(\bar{X}(t), \bar{u}(t)) Q_j(t)+\sum_{j=1}^d Q_j(t) \sigma_x^j(\bar{X}(t), \bar{u}(t))+H_{x x}(\bar{X}(t), \bar{u}(t), p(t), K(t))\right] \mathrm{d}t\\
&-Q(t) \mathrm{d} W(t),\\
		 P(\tau^{\bar{u}})=&g_{x x}(\bar{X}(\tau^{\bar{u}})). \\
	\end{aligned}
\end{equation}

By variational equations (\ref{eq-1}), (\ref{eq-2}) and first-order adjoint equation (\ref{first-or}), we can rewrite (\ref{vi-4}) as
 \begin{equation}
 	\label{eq-5}
 	 \begin{aligned}
 	 	&\mathbb{E} \int_{0}^{\tau^{\bar{u}}}\left(H\left(\bar{X}(s), u^{\varepsilon}(s), p(s), K(s)\right)-H(\bar{X}(s), \bar{u}(s), p(s), K(s))\right) \mathrm{d}s+\mathbb{E}\left[\varepsilon\frac{k(\tau)}{h^{\bar{u}}(\tau^{\bar{u}})}R(\tau^{\bar{u}})\right]\\
 		& +\frac{1}{2} \mathbb{E} \int_{0}^{\tau^{\bar{u}}} y_1^*(s) H_{x x}(\bar{X}(s), \bar{u}(s), p(s), K(s)) y_1(s)\mathrm{d}s+\frac{1}{2} \mathbb{E}[ y_1^*(\tau^{\bar{u}})g_{x x}(\bar{X}(\tau^{\bar{u}})) y_1(\tau^{\bar{u}})]  \geq o(\varepsilon),
 	\end{aligned}
 \end{equation}
 where
 $$
 H(x, v, p, k)=f(x, v)+(p, b(x, v))+\sum_{j=1}^d\left(k_j, \sigma^j(x, v)\right).
 $$

Based on the second-order adjoint equation (\ref{second-or}), we further rewrite (\ref{eq-5}) as
%
 $$
 \begin{aligned}
 	&\mathbb{E}\int_{0}^{\tau^{\bar{u}}}\left(H\left(\bar{X}(t), u^{\varepsilon}(t), p(t), K(t)\right)-H(\bar{X}(t), \bar{u}(t), p(t), K(t))\right) \mathrm{d}t
 +\mathbb{E}\left[\varepsilon\frac{k(\tau)}{h^{\bar{u}}(\tau^{\bar{u}})}R(\tau^{\bar{u}})\right] \\
 	&+\frac{1}{2} \mathbb{E} \int_{0}^{\tau^{\bar{u}}} \operatorname{tr}\left[\left(\sigma\left(\bar{X}(t), u^{\varepsilon}(t)\right)-\sigma(\bar{X}(t), \bar{u}(t))\right)^* P(s)\left(\sigma\left(\bar{X}(t), u^{\varepsilon}(t)\right)-\sigma(\bar{X}(t),\bar{u}(t))\right)\right] \mathrm{d}t \geq o(\varepsilon).\\
 \end{aligned}
 $$
 Thus, we obtain the following variational inequality
 $$
 \begin{aligned}
 	& H(\bar{X}(\tau), u, p(\tau), K(\tau))-H(\bar{X}(\tau), \bar{u}(\tau), p(\tau), K(\tau))
 +\frac{k(\tau)}{h^{\bar{u}}(\tau^{\bar{u}})}R(\tau^{\bar{u}})\\
 	&+\frac{1}{2} \operatorname{tr}\left[(\sigma(\bar{X}(\tau), u)-\sigma(\bar{X}(\tau), \bar{u}(\tau)))^* P(\tau)(\sigma(\bar{X}(\tau), u)-\sigma(\bar{X}(\tau), \bar{u}(\tau)))\right]
 	\geq  0 \\
 & \forall u \in U, \quad \text { a.e., a.s. }
 \end{aligned}
 $$

 The cases (ii) and (iii) of Lemma \ref{le-5} are similar to the case (i). Thus, we omit the details of cases (ii) and (iii). We conclude the main result of this paper as follows.
 \begin{theorem}\label{th-1}
 Let Assumptions \ref{ass-a} and \ref{ass-b} hold, and suppose that $h^{\bar{u}}(\tau^{\bar{u}})\neq 0$ and $h^{\bar{u}}(\cdot)$ is continuous at the point $\tau^{\bar{u}}$. If $(\bar{X}(\cdot), \bar{u}(\cdot))$ is an optimal pair of this control system.  Then we have
 $$
 \begin{aligned}
 	& (p(\cdot), K(\cdot)) \in L_{\mathcal{F}}^{2}(0,\tau^{\bar{u}};\mathbb{R}^{m})\times(L_{\mathcal{F}}^{2}(0,\tau^{\bar{u}};\mathbb{R}^{m}))^{d}, \\
 	& (P(\cdot), Q(\cdot)) \in  L_{\mathcal{F}}^2\left(0, \tau^{\bar{u}} ; \mathbb{R}^{m\times m}\right) \times\left(L_{\mathcal{F}}^2\left(0, \tau^{\bar{u}} ; \mathbb{R}^{m\times m}\right)\right)^d, \\
 \end{aligned}
 $$
 and for $0\leq \tau\leq \tau^{\bar{u}}$ the following results hold:

 (i). If $\tau^{\bar{u}}<T$, one obtains
  \begin{equation}
   \begin{aligned}
  	& H(\bar{X}(\tau), u, p(\tau), K(\tau))-H(\bar{X}(\tau), \bar{u}(\tau), p(\tau), K(\tau))
  +\frac{k(\tau)}{h^{\bar{u}}(\tau^{\bar{u}})}R(\tau^{\bar{u}})\\
  &+\frac{1}{2} \operatorname{tr}\left[(\sigma(\bar{X}(\tau), u)-\sigma(\bar{X}(\tau), \bar{u}(\tau)))^* P(\tau)(\sigma(\bar{X}(\tau), u)-\sigma(\bar{X}(\tau), \bar{u}(\tau)))\right]
  \geq  0 \\
  & \forall u \in U, \ \text { a.e., a.s. }
  \end{aligned}	
  \end{equation}

 (ii). If $\inf\bigg{\{}t:\mathbb{E}[\Phi(\bar{X}(t))]\geq \alpha,\ t\in [0,T] \bigg{\}}=T$, one obtains

 \begin{equation}
 	\begin{aligned}
 		& H(\bar{X}(\tau), u, p(\tau), K(\tau))-H(\bar{X}(\tau), \bar{u}(\tau), p(\tau), K(\tau)) +\frac{k(\tau)}{h^{\bar{u}}(\tau^{\bar{u}})}R(\tau^{\bar{u}})\\
 	&+\frac{1}{2} \operatorname{tr}\left[(\sigma(\bar{X}(\tau), u)-\sigma(\bar{X}(\tau), \bar{u}(\tau)))^* P(\tau)(\sigma(\bar{X}(\tau), u)-\sigma(\bar{X}(\tau), \bar{u}(\tau)))\right]
 	\geq  0 \\
 & \forall u \in U, \ \text { a.e., a.s. }
 	\end{aligned}	
 \end{equation}

or
\begin{equation}
	\begin{aligned}
		& H(\bar{X}(\tau), u, p(\tau), K(\tau))-H(\bar{X}(\tau), \bar{u}(\tau), p(\tau), K(\tau)) \\
		&+\frac{1}{2} \operatorname{tr}\left[(\sigma(\bar{X}(\tau), v)-\sigma(\bar{X}(\tau), \bar{u}(\tau)))^* P(\tau)(\sigma(\bar{X}(\tau), u)-\sigma(\bar{X}(\tau), \bar{u}(\tau)))\right]
		 \geq  0 \\
& \quad \forall u \in U, \ \text { a.e., a.s. }
	\end{aligned}	
\end{equation}

 (iii). If $\bigg{\{}t:\mathbb{E}[\Phi(\bar{X}(t))]\geq \alpha,\ t\in [0,T] \bigg{\}}=\varnothing$, we have
 \begin{equation}
 	\begin{aligned}
 		& H(\bar{X}(\tau), u, p(\tau), K(\tau))-H(\bar{X}(\tau), \bar{u}(\tau), p(\tau), K(\tau)) \\
 		&+\frac{1}{2} \operatorname{tr}\left[(\sigma(\bar{X}(\tau), u)-\sigma(\bar{X}(\tau), \bar{u}(\tau)))^* P(\tau)(\sigma(\bar{X}(\tau), v)-\sigma(\bar{X}(\tau), \bar{u}(\tau)))\right]
 		 \geq  0\\
 & \forall u \in U, \ \text { a.e., a.s. }
 	\end{aligned}	
 \end{equation}
\end{theorem}

In the following, we give two examples to verify our main result. The first one is about deterministic optimal system with varying terminal time and the second one is about stochastic system.
\begin{example}\label{exm-1}
For a two-point control set $U=\{1,2\}$, let $m=d=1,T=1,\ \alpha=1$, $b(x,u)=x+u,\ \sigma(x,u)=0,\ f(x,u)=u, \ \Phi(x)=x$ and $g(x)=0$. The state equation is given as follows:
\begin{equation}
	X^u(s)=\int_0^s\bigg{[}X^u(t)+u(t)\bigg{]}\mathrm{d}t,
\end{equation}
where the varying terminal time is
$$
\tau^{{u}}=\inf\bigg{\{}t:{X}^u(t)\geq \alpha,\ t\in [0,T] \bigg{\}}\bigwedge T,
$$
and the cost functional is
$$
J(u(\cdot))=\int_0^{\tau^{{u}}}{u}(s)\mathrm{d}s.
$$
Then for $t\in[0,T]$, we can find an optimal pair $(\bar{u}(t),\bar{X}(t))$ by using Theorem $\ref{th-1}$. We suppose that $\tau^{\bar{u}}<1$. For $t\in[0,\tau^{\bar{u}}]$, we have $h^{\bar{u}}(t)=\bar{X}(t)+\bar{u}(t)$, and
first- and second- order variational equations
$$
y_{1}(t)=\int_{0}^{t}\bar{u}(s)y_{1}(s)+(u^{\varepsilon}(s)-\bar{u}(s))ds,
$$
$$
y_{2}(t)=\int_{0}^{t}y_{2}(s)+(u^{\varepsilon}(s)-\bar{u}(s))ds.
$$
The corresponding first- and second- order adjoint equations are
$$
 \begin{aligned}
 &-\mathrm{d}p(t)=-K(t)\mathrm{d}W(t)\\
 &p(\tau^{\bar{u}})=0,\\	
 \end{aligned}
$$
$$
 \begin{aligned}
 &-\mathrm{d}P(t)=-Q(t)\mathrm{d}W(t)\\
 &P(\tau^{\bar{u}})=0\\	
 \end{aligned}
$$
and we have
$$
\begin{aligned}
&-\mathrm{d}p_{0}(t)=-K_{0}(t)\mathrm{d}W(t)\\
&p_{0}(\tau^{\bar{u}})=0,\\	
\end{aligned}
$$
$$
\begin{aligned}
&-\mathrm{d}P_{0}(t)=-Q_{0}(t)\mathrm{d}W(t)\\
&P_{0}(\tau^{\bar{u}})=0.\\	
\end{aligned}
$$
Thus,we can obtain
\begin{equation}\label{eq-*}
k(\tau)=u-\bar{u}(\tau)	
\end{equation}
By (i) of Theorem $\ref{th-1}$, combining $(\ref{eq-*})$ and $\bar{X}(\tau^{\bar{u}})=\alpha=1$, we have
$$
\begin{aligned}
	& H(\bar{X}(\tau), u, p(\tau), K(\tau))-H(\bar{X}(\tau), \bar{u}(\tau), p(\tau), K(\tau))
	+\frac{k(\tau)}{h^{\bar{u}}(\tau^{\bar{u}})}R(\tau^{\bar{u}})\\
	&+\frac{1}{2} \operatorname{tr}\left[(\sigma(\bar{X}(\tau), u)-\sigma(\bar{X}(\tau), \bar{u}(\tau)))^* P(\tau)(\sigma(\bar{X}(\tau), u)-\sigma(\bar{X}(\tau), \bar{u}(\tau)))\right]\\
	&=(u-\bar{u}(\tau))\big[1+\frac{\bar{u}(\tau^{\bar{u}})}{\bar{X}(\tau^{\bar{u}})+\bar{u}(\tau^{\bar{u}})}\big]\geq  0. \\
\end{aligned}
$$

Suppose that $1+\frac{\bar{u}(\tau^{\bar{u}})}{\bar{X}(\tau^{\bar{u}})+\bar{u}(\tau^{\bar{u}})}\geq  0$, we can obtain an optimal pair $\big(\bar{u}(t),\bar{X}(t)\big)=(1,e^{t}-1),\forall t\in[0,\tau^{\bar{u}}]$.

According to $\bar{X}(\tau^{\bar{u}})=e^{\tau^{\bar{u}}}-1=1$, we have $\tau^{\bar{u}}=\ln2<T=1$. Notice that $1+\frac{\bar{u}(\tau^{\bar{u}})}{\bar{X}(\tau^{\bar{u}})+\bar{u}(\tau^{\bar{u}})}=\frac{3}{2}\geq 0$, thus we find the optimal pair $\big( \bar{u}(t),\bar{X}(t)\big)=(1,e^{t}-1),\forall t\in[0,\tau^{\bar{u}}].$

Similarly, for a convex control set $U=[1,2]$ according to Theorem $\ref{th-1}$, we still can obtain an optimal pair $\big(\bar{u}(t),\bar{X}(t)\big)=(1,e^{t}-1),\forall t\in[0,\tau^{\bar{u}}]$.\\
\end{example}
\begin{example}
Let $U=\{1,2\}$, $m=d=1,T=1,\ \alpha=1$, $b(x,u)=1,\ \sigma(x,u)=u,\ f(x,u)=u, \ \Phi(x)=x$ and $g(x)=0$. The state equation is given as follows:
\begin{equation}
	X^u(t)=\frac{1}{2}+t+\int_0^t u(s)\mathrm{d}W(s),
\end{equation}
where the varying terminal time is
$$
\tau^{{u}}=\inf\bigg{\{}t:\mathbb{E}[{X}^u(t)]\geq 1,\ t\in [0,T] \bigg{\}}\bigwedge T,
$$
and the cost functional is
$$
J(u(\cdot))=\int_0^{\tau^{{u}}}{u}(s)\mathrm{d}s.
$$
In this case, the first- and second-order variantial equations are
$$
\begin{aligned}
	&y_{1} (t)=\int_{0}^{t} u^{\varepsilon}(s)-\bar{u}(s) \mathrm{d}s,\\
	&y_{2} (t)=0.\\
	\end{aligned}
$$
The adjoint equations are same with those given in Example $\ref{exm-1}$. Therefore we have $l(X^u(t),u(t))=1$ and $k(t)=0,\ \forall t\in[0,\tau^{\bar{u}}]$.\\

Suppose that $\tau^{\bar{u}}<1$, by Theorem $\ref{th-1}$, we have
$$
\begin{aligned}
	& H(\bar{X}(\tau), u, p(\tau), K(\tau))-H(\bar{X}(\tau), \bar{u}(\tau), p(\tau), K(\tau))
	+\frac{k(\tau)}{h^{\bar{u}}(\tau^{\bar{u}})}R(\tau^{\bar{u}})\\
	&+\frac{1}{2} \operatorname{tr}\left[(\sigma(\bar{X}(\tau), u)-\sigma(\bar{X}(\tau), \bar{u}(\tau)))^* P(\tau)(\sigma(\bar{X}(\tau), u)-\sigma(\bar{X}(\tau), \bar{u}(\tau)))\right]\\
	&=u-\bar{u}(\tau)\geq  0. \\
\end{aligned}
$$
Thus we can obtain an optimal control $\bar{u}(\tau)=1$ and optimal state $\bar{X}(\tau)=\frac{1}{2}+\tau+W(\tau),\ \forall\tau \in[0,\ \tau^{\bar{u}}]$,
which deduces that $\tau^{\bar{u}}=\frac{1}{2}<1$.
\end{example}

\section{Conclusion}
In this study, we consider a varying terminal time structure for the stochastic optimal control problem under state constraints, in which the terminal time varies with the mean value of the state via the constrained condition. When the control domain is non-convex, the related global stochastic maximum principle is still unsolved. We establish the related global stochastic maximum principle for the optimal control problem with varying terminal time. Following the idea given in Peng \cite{P90} and Yang \cite{Y20}, we use the spike variation method to develop the stochastic maximum principle.  To overcome the difficulty in the optimal control problem, we establish the variational equation for the varying terminal time and cost functional based on three kinds of case of optimal terminal time $\tau^{\bar{u}}$: $\text{(i)} \ \tau^{\bar{u}}<T$;
$\text{(ii)}\ \inf\bigg{\{}t:\mathbb{E}[\Phi(\bar{X}(t))]\geq \alpha,\ t\in [0,T] \bigg{\}}=T$;
$\text{(iii)}\ \bigg{\{}t:\mathbb{E}[\Phi(\bar{X}(t))]\geq \alpha,\ t\in [0,T] \bigg{\}}=\varnothing$. Based on the developed variational equation of the varying terminal time, first- and second-order adjoint equations, we establish a global stochastic maximum principle for our optimal control problem.

\section*{Statements \& Declarations}

\textbf{Funding}: This work was supported by the National Key R\&D program of China (Grant No.2023YFA1009203, 2018YFA0703900), National Natural Science Foundation of China (Grant No.11701330), and Taishan Scholar Talent Project Youth Project.

\noindent \textbf{Competing Interests}: The authors have no relevant financial or non-financial interests to disclose.

\noindent\textbf{Author Contributions}: Shuzhen Yang contributed to the study conception and design. The first draft of the manuscript was written by Jin Shi and all authors commented on previous versions of the manuscript. All authors read and approved the final manuscript.
\end{document}